\documentclass[11pt]{article}
\usepackage{amsfonts}
\usepackage{color}
\usepackage{graphics}

\usepackage{cite}
\usepackage{latexsym}
\usepackage{amsmath}
\usepackage{amssymb}
\usepackage[dvips]{epsfig}
\usepackage{amscd}

 \usepackage{color}
 \usepackage{indentfirst}
\usepackage{ntheorem}
\theoremstyle{plain}
\theoremseparator{.}
\newtheorem{theorem}{Theorem}[section]
\newtheorem{remark}{Remark}[section]

\newtheorem{lemma}[theorem]{Lemma}

\newcommand\thmref[1]{Theorem~\ref{#1}}
\newcommand\lemref[1]{Lemma~\ref{#1}}

\newcommand{\nb}{{\bold b}}
\newcommand{\nw}{{\bold w}}
\newcommand{\ti}{\tilde}

\newcommand{\lm}{\lambda}

\def\pf{{\it Proof.}  }

\newcommand{\thatsall}{\hfill$\Box$}
\newcommand{\bi}{\bibitem}

\newcommand{\bt}{\begin{theorem}}
\newcommand{\bl}{\begin{lemma}}
\newcommand{\el}{\end{lemma}}
\newcommand{\et}{\end{theorem}}

\renewcommand{\b}{\beta  }

\newcommand{\te}{\theta}

\newcommand{\al}{\alpha}

\newcommand{\la}{\label}

\newcommand{\ka}{\kappa}

\newcommand{\bn}{\begin{eqnarray}}
\newcommand{\en}{\end{eqnarray}}
\newcommand{\bnn}{\begin{eqnarray*}}
\newcommand{\enn}{\end{eqnarray*}}

\newcommand{\bnnn}{\begin{eqnarray*}}
\newcommand{\ennn}{\end{eqnarray*}}
\newcommand{\ben}{\begin{enumerate}}
\newcommand{\een}{\end{enumerate}}
\newcommand{\ba}{\begin{aligned}}
\newcommand{\ea}{\end{aligned}}
\newcommand{\be}{\begin{equation}}
\newcommand{\ee}{\end{equation}}

\def\norm[#1]#2{\|#2\|_{#1}}

\def\O{\Omega}

\newcommand{\xl}{\left}
\newcommand{\xr}{\right}

\makeatletter
\@addtoreset{equation}{section}
\makeatother

\makeatletter
\@addtoreset{equation}{section}
\makeatother

\title{ Global existence of strong solutions to the planar compressible magnetohydrodynamic equations with large initial data in unbounded domains\thanks{ \textsc{B. L\"u} is supported by NNSFC (11971217) and Jiangxi Provincial Natural Science Foundation (20202ACBL211002). \textsc{X. Shi} is supported by NNSFC (11671027 \&11471321).
}
}

\author{ Boqiang L\"u\thanks{College of Mathematics and Information Science, Nanchang Hangkong University, Nanchang 330063, P. R. China ({\tt lvbq86@163.com}). }
\quad Xiaoding Shi \thanks{Department of Mathematics, College of Mathematics and Physics, Beijing University of Chemical Technology, Beijing 100029, P.R. China ({\tt shixd@mail.buct.edu.cn}).
 }
\quad Chengfeng  Xiong\thanks{Institute of Applied Mathematics, AMSS,
Chinese Academy of Sciences, Beijing 100190,  P. R. China ({\tt xiongcf1998@163.com}).
 }
 }
\date{ }

\begin{document}
\maketitle

\begin{abstract}
In one-dimensional unbounded domains, we consider the equations of a planar compressible magnetohydrodynamic (MHD) flow with constant viscosity and heat conductivity. More precisely,  we prove the global existence of strong solutions to the MHD equations with large initial data satisfying the same conditions as those of Kazhikhov's theory  in bounded domains (Kazhikhov 1987 Boundary Value Problems for Equations of Mathematical Physics (Krasnoyarsk)). In particular, our result generalizes the Kazhikhov's theory for the initial boundary value problem in bounded domains to the unbounded case.
  \end{abstract}

$\mathbf{Keywords.}$ Magnetohydrodynamics; Global strong solutions; Large initial data; Unbounded domains

{\bf Math Subject Classification:} 35Q35; 76N10.

\section{Introduction}
Magnetohydrodynamics (MHD) is concerned with the study of the interaction between the magnetic fields and electrically conducting fluids. It is widely applied to astrophysics, geophysics and plasma physics in practice,
see \cite{q1,q2,q3,q4,q5,q6,q7}
and the references therein.  We are concerned with the governing equations of a planar magnetohydrodynamic compressible flow written in the Lagrange variables

\be\la{1.1}
v_t=u_{x},
\ee
\be\la{1.2}
u_{t}+(P+\frac12 |\nb|^2)_{x}=\left(\mu\frac{u_{x}}{v}\right)_{x},\ee
 \be\la{1.3}\nw_t-\nb_x=\left(\lm\frac{\nw_{x}}{v}\right)_{x}, \ee \be\la{1.4} (v\nb)_t-\nw_x =\left(\nu\frac{\nb_{x}}{v}\right)_{x}, \ee
\be\la{1.5}\ba&\left(e+\frac{u^2+|\nw|^2+v|\nb|^2}{2}\right)_{t}+ \left(
u\left(P+\frac12|\nb|^2\right)-\nw\cdot\nb\right)_{x}\\&\quad=\left(\kappa\frac{\theta_{x}}{v}
+\mu\frac{uu_x}{v}+\lm\frac{\nw\cdot\nw_x}{v}+\nu\frac{\nb\cdot\nb_x}{v}\right)_x,
 \ea
\ee
where $t>0$ is time, $x\in \Omega \subset  \mathbb{R}=(-\infty,+\infty)$ denotes the
Lagrange mass coordinate,  and the unknown functions $v>0,~u,~\nw\in \mathbb{R}^2,~\nb \in   \mathbb{R}^2, e>0, \theta>0,$ and $P$ are,  respectively, the specific volume of the gas, longitudinal velocity, transverse velocity, transverse magnetic field, internal energy,  absolute temperature and  pressure. $\mu$ and $\lm$ are the viscosity of the flow,
$\nu$ is the magnetic diffusivity of the magnetic field, and $\ka$ is the heat conductivity.
In this paper, we
consider a perfect gas for magnetohydrodynamic flow, that is, $P$ and $e$ satisfy \be \la{1.6}   P =R \theta/{v},\quad e=c_v\theta +\mbox{const},
\ee
where  both specific gas constant  $R$ and   heat
capacity at constant volume $c_v $ are   positive constants.
We also assume that $\lm$ and  $\nu$  are  positive constants, and that $\mu,~\ka$ satisfy
\be\label{p1.8} \mu=\ti\mu_1+\ti\mu_2v^{-\alpha},~~~~\ka=\ti\ka\te^{\beta},\ee
with constants $\ti\mu_1>0,~\ti\mu_2\ge 0,~\ti\ka>0$, and $\alpha,~\beta\ge0$.

The system \eqref{1.1}-\eqref{1.6} is supplemented
with the initial  conditions
\be\la{1.8}(v,u,\te,\nb,\nw)(x,0)=(v_0,u_0,\te_0,\nb_0,\nw_0)(x),  \quad  x\in\Omega, \ee
and one of three types of far-field and boundary  ones:

 1) Cauchy problem
\be\la{1.9}\Omega=\mathbb{R}, \quad\lim_{|x|\rightarrow\infty}\left(v,u,\theta,\nb,\nw\right)=(1,0,1,0,0),\quad t>0;\ee

 2) boundary and far-field conditions for $\Omega=(0,\infty)$
\be\ba\la{1.10} &u(0,t)=0,\quad\theta(0,t)=1,\quad\nb(0,t)=\nw(0,t)=0,\,\\&\lim_{x\rightarrow\infty}\left(v,u,\theta,\nb,\nw\right)=(1,0,1,0,0),\quad t>0; \ea\ee

 3) boundary and far-field conditions for $\Omega=(0,\infty)$
\be\ba\la{1.11} &u(0,t)=\theta_x(0,t)=0,\quad\nb(0,t)=\nw(0,t)=0,\\& \lim_{x\rightarrow\infty}\left(v,u,\theta,\nb,\nw\right)=(1,0,1,0,0),\quad t>0. \ea\ee


There is huge literature on the studies of the global existence and large time behavior of solutions to the compressible MHD system. In particular, when $\nw = \nb =  0$, the   MHD system \eqref{1.1}-\eqref{1.5} reduces to the Navier-Stokes equations, which have been studied extensively in \cite{L-L,9,6,cw1,cw2,fjn1,25,fhl1,kazh,py,24,28,hs1} and the references therein.  Here, we recall briefly  some results  which
are more relative with our problem.   For constant coefficients $\al=\beta=0$ and large initial data, Kazhikhov-Shelukhin \cite{9} first obtained the global existence of solutions to the initial boundary value problem in bounded domains. When $\al=0$ and $\b>0$ in \eqref{p1.8},   Huang-Shi \cite{hs1} obtained the global strong solutions to the initial-boundary-value problem with the initial data  $(v_0,u_0,\te_0)\in  H^1$, see also  Pan-Zhang \cite{28} for the initial data $(v_0,u_0,\te_0)\in (H^1\times H^2\times H^2)$. For the Cauchy problem in the unbounded domain,  Kazhikhov \cite{kazh}   obtained the global existence of strong solutions  with  constant coefficients $\al=\beta=0$, which recently be refined  by Li-Shu-Xu \cite{py}  to the general case with $\al=0$ and $\b>0$.

 Now, let's  go back to the MHD system   \eqref{1.1}-\eqref{1.5}. For the  initial-boundary-value problem in bounded domains, the global existence of strong solutions with large initial data was obtained by  Kazhikhov  \cite{ka1} (see also Amosov-Zlotnik  \cite{az1}) for  the  constant coefficients $\al=\beta=0$ and by Huang-Shi-Sun \cite{hss1} for $\al\ge0$ and $\b>0$ (see also Hu-Ju \cite{7} with $\al=0$ and $\b>0$).
 Concerning the unbounded domains,  Cao-Peng-Sun \cite{cps} established the global existence of strong solutions with large initial data.
 It should be pointed that the method in \cite{cps} depends heavily on the
 assumption of $\beta >0$  and thus cannot be
 adapted  to the case of constant coefficients $\al=\beta=0$.
  Therefore, the main aim of this paper is to prove the global existence of strong solutions with     constant coefficients $\al=\beta=0$ in unbounded domains, which
  generalized  Kazhikhov's result  \cite{ka1} to the  case of unbounded domains.   That is, our main result is as follows.

  \begin{theorem}\la{thm1.1} Suppose that $\alpha=\beta=0$ and the initial data $ ( v_0,u_0,\te_0,\nb_0,\nw_0)$   satisfies
  \be  \la{a1.11} ( v_0-1,u_0,\te_0-1,\nb_0,\nw_0)\in H^1(\O), \ee  and \be \la{a1.12}
\inf_{x\in \O}v_0(x)>0, \quad \inf_{x\in \O }\theta_0(x)>0, \ee
and are compatible with \eqref{1.10}, \eqref{1.11}. Then there exists a unique global strong solution $(v,u,\theta,\nb,\nw)$  to the initial-boundary-value problem \eqref{1.1}-\eqref{1.9}, or \eqref{1.1}-\eqref{1.8} \eqref{1.10}, or \eqref{1.1}-\eqref{1.8} \eqref{1.11} satisfying for any $T>0$,
 \be
 \begin{cases} \la{a1.13}  v-1, u,\,\theta-1,\nb,\nw \in L^\infty(0,T;H^1(\O)),\\ v_t\in
  L^\infty(0,T;L^2(\O))\cap L^2(0,T;H^1(\O)), \\ u_t,\,\theta_t,\,\nb_t,\,\nw_t,\,u_{xx},\,\te_{xx},\,\nb_{xx},\,\nw_{xx} \in
  L^2((\O)\times(0,T)),\end{cases}\ee
  and for each $(x,t)\in \O \times[0,T]$
  \be C^{-1}\leq v(x,t)\leq C,\quad C^{-1}\leq\te(x,t)\leq C,\ee
  where $C>0$ is a constant  depending on the   data and T.
  \end{theorem}

\begin{remark} Our result can be regarded as a natural generalization of    Kazhikhov's result  \cite{ka1} in the bounded domains to the  case of unbounded domains.  \end{remark}

 \begin{remark}  It should be mentioned here that Theorem \ref{thm1.1} still holds for the case of $\nw = \nb = 0$. This in particular yields that our result also establishes  the global strong solutions for compressible Navier-Stokes equations which has been considered   in   Kazhikhov \cite{kazh}.   \end{remark}

We now comment on the analysis of this paper.
To extend the local strong solutions whose existence is guaranteed  by lemma \ref{plq1} to be global,  the main issue is to establish some necessary global a priori estimates of solutions.
 Motivated by Kazhikhov \cite{ka1} (see also   \cite{hss1,hs1}), we first  obtain a key representation of $v$ (see \eqref{2.000}), which together with the standard energetic estimates \eqref{2.2}  derives the lower  bound  of $v$   \eqref{cqq1}. Then, following the similar arguments as those in \cite{py},   multiplying
the temperature equation \eqref{2.1} by $\te^{-2}(\te^{-1}-2)^p_+$ and  using the boundedness of domains  $(\te<1/2)(t)$ (see  \eqref{2.18}),
     we can also obtain   the lower bound of $\te$ (see \eqref{cqq1x}). Next, we will prove  the key   upper bound  of  $v$. It should be mentioned here that the method in \cite{cps} for the case $\beta>0$ is not  valid  in the case $\beta=0$ due to lack of the estimates on $L^1(0,T; L^\infty(\O))$-norm of $\te$. In this paper,
      modifying slightly the idea due to  Amosov-Zlotnik \cite{az1},  we can  prove the  key upper bound  of  $v$ (see \eqref{q2.49})  by controlling the  $L^\infty(0,T;L^2)$-norm of $(\ln v)_x$ (see \eqref{xiao6}), see Lemma \ref{upv}. Finally, using  the similar arguments as those in \cite{py,cps,hss1},   one can  derive the necessary a priori estimates of the solutions, see Lemmas \ref{lemmy0}--\ref{nlemma80}.   The whole procedure will be carried out in the next section.

\section{ Proof of \thmref{thm1.1}}

We first state   the following  existence and uniqueness of local solutions which can be obtained by using the Banach theorem and the contractivity of the operator defined by the linearization of the problem on a small time interval (c.f. \cite{10,13,tan}).

\begin{lemma} \la{plq1} Under the assumptions in Theorem \ref{thm1.1},  there exists some $T_0>0$ such that  the initial-boundary-value problem \eqref{1.1}-\eqref{1.9}, or \eqref{1.1}-\eqref{1.8} \eqref{1.10}, or \eqref{1.1}-\eqref{1.8} \eqref{1.11} has a unique strong solution $(v,u,\te,\nb,\nw)$ with positive $v(x,t)$ and $\theta(x,t)$ satisfying \eqref{a1.13}.
\end{lemma}

Then, to finish the proof of  Theorem \ref{thm1.1},
it only remains to obtain some a priori estimates (see \eqref{cqq1}, \eqref{cqq1x}, \eqref{q2.49}, \eqref{uiz1}, \eqref{aas1}, \eqref{7.1}, and \eqref{eq1}),
where the constants depend only on $T$ and the data of the problem. Thus, one can use the a priori estimates to continue the local solutions to the whole interval $[0, T]$.

Next, without loss of generality, we assume that $\lambda=\nu= \mu = \ti\ka=R=c_v=1$.
We first state  the following  basic energy estimates.
\begin{lemma}\la{enest} It holds that for any $(x,t)\in \Omega\times [0,T],$
\be\ba\la{2.2}&\sup_{0\le t\le T}\int_\O\left(\frac{u^2+|\nw|^2+v|\nb|^2}{2}+(v-\ln v-1) +(\theta-\ln \theta-1)\right)dx
 +\int_0^T W(t)dt\le e_0,      \ea\ee
where \be\la{2.3}W(t)   \triangleq \int_\O\left(\frac{ \theta_{x}^2}{v\theta^2 }+\frac{u_x^2+|\nw_x|^2+|\nb_x|^2}{v\theta}\right)dx, \nonumber\ee
and
\be e_{0}\triangleq  2\int_\O\left(\frac{ u_0^2+|\nw_{0}|^2+v_0|\nb_{0}|^2}{2}+(v_0-\ln v_0-1)+(\theta_0-\ln \theta_0-1)\right)dx\nonumber . \ee
\end{lemma}

\pf Using  \eqref{1.1}-\eqref{1.4},  the energy equation \eqref{1.5}  can be  rewritten as
\be\la{2.1}\theta_{t}+  \frac{\theta}{v}u_{x}= \left(\frac{ \theta_{x}}{v}\right)_{x}+ \frac{  u_{x}^{2}+|\nw_x|^2+|\nb_x|^2}{v}.\ee
Multiplying   \eqref{1.1} by $ 1- {v}^{-1}$, \eqref{1.2} by $u$, \eqref{1.3} by $\nw$, \eqref{1.4} by $\nb$, and \eqref{2.1} by $  1- {\theta}^{-1} $,  respectively, one obtains  after  adding the resultant equalities altogether that
\bnn\ba &\left(\frac{u^2+|\nw|^2+v|\nb|^2}{2}+\left(\theta-\ln \theta-1\right)+\left(v-\ln v-1\right)\right)_{t}\\&\quad +\frac{  \theta_{x}^2}{v\theta^2}+\frac{ u_x^2+|\nw_{x}|^2+|\nb_{x}|^2}{v\theta}\\
&=\left(\frac{ \theta_{x}}{v}+\frac{u u_x}{v}+\frac{\nw\cdot\nw_{x}}{v}+\frac{\nb\cdot\nb_{x}}{v}\right)_{x}+u_x \\&\quad-\left(u\left(\frac{\theta}{v}+\frac{1}{2}|\nb|^2\right)
-\nw\cdot\nb\right)_{x}-\left(\frac{ \theta_{x}}{v\theta}\right)_x,
\ea\enn
which along with \eqref{1.9} or \eqref{1.10} or \eqref{1.11} yields \eqref{2.2}and completes the proof Lemma \ref{enest}. \thatsall

Next, we will give out the following   lower bound   of $v$.
\begin{lemma}\la{wq20} There  exists a positive constant $C$ such  that  for any $(x,t)\in\O\times [0,T]$,
\be\ba\la{cqq1}
v(x,t)\ge C,
\ea\ee
where (and in what follows)   $C $  denotes a
	generic positive constant
	depending only on $T,  \|(v_0-1,u_0,\theta_0-1,  \nb_0,\nw_0)\|_{H^1(\O)},
	\inf\limits_{x\in \O}v_0(x),$ and $ \inf\limits_{x\in \O}\theta_0(x).$
\end{lemma}

\pf Letting
\bnn\ba\la{2.9}
\sigma&\buildrel\Delta\over=\frac{ u_x}{v}-\left(\frac{\theta}{v}+\frac{1}{2}|\nb|^2\right) =(\ln v)_t-\left(\frac{\theta}{v}+\frac{1}{2}|\nb|^2\right)
\ea\enn
owing to \eqref{1.1},
we write \eqref{1.2} as
\be\ba\la{2.8} u_t=\sigma_x. \ea\ee

For any $x\in\O$, denoting $N=[x],$ one obtains after
integrating \eqref{2.8} over $[N,x]\times[0,t]$ that
\be\ba
\int_N^{x}u dy-\int_N^{x}u_0 dy=&\ln v-\ln v_0\\&-\int_0^t\left(\frac{\theta}{v}+\frac{1}{2}|\nb|^2 \right)d\tau
-\int_0^t\sigma(N,\tau)d\tau,\nonumber
\ea\ee
 which implies
\be\la{2.11}\ba v(x,t) =   B_N(x,t)Y_N(t)  \exp\left\{ \int_0^t\left(\frac{\te}{v}+\frac{1}{2}|\nb|^2\right)d\tau\right\}, \ea\ee
where
\be\ba \la{2.0}B_N(x,t)\triangleq v_0 \exp\left\{\int_N^x u dy-\int_N^x u_0 dy\right\},  \ea\ee
and
\be\ba\la{2.00} Y_N(t)\triangleq\exp\left\{\int_0^t \sigma(N,\tau)d\tau\right\}. \ea\ee
 Denoting by
 \be\ba g(x,t)=\int_0^t\left(\frac{\te+\frac{1}{2}v|\nb|^2}{v}\right)d\tau,\nonumber\ea\ee
it deduces from \eqref{2.11} that
 \be\ba g_t=\frac{\theta+\frac{1}{2}v|\nb|^2}{v}=\frac{\theta+\frac{1}{2}v|\nb|^2}{B_N(x,t)Y_N(t)\exp\left\{g\right\}},\nonumber\ea\ee
 which gives
 \be\ba \exp\left\{g\right\}=1+\int_0^t\frac{\te+\frac{1}{2}v|\nb|^2}{B_N(x,\tau)Y_N(\tau)}d\tau .\nonumber \ea\ee
 Combining this with  \eqref{2.11} leads to
\be \la{2.000}\ba v(x,t)=   B_N (x,t)Y_N (t)\left(1+\int_0^t\frac{(\te+\frac{1}{2}v|\nb|^2)}{B_N(x,\tau)Y_N(\tau)}d\tau \right).\ea\ee

On the one hand, it follows from \eqref{2.2}  that
\be\ba\left|\int_N^x\left(u(y,t)-u_0(y)\right)dy\right|\le\left(\int_N^{N+1} u^2 dy\right)^{\frac{1}{2}}+\left(\int_N^{N+1}u_0^2dy\right)^{\frac{1}{2}}\le C,\nonumber\ea\ee
which together with \eqref{2.0} implies
\be\la{2.13}    C^{-1}\le B_N(x,t)\le C,  \ee
where, and in what follows, $C$ is a constant independent of $N$.

On the other hand,  it follows from \eqref{2.2} that
$$ \int_{N}^{N+1} (v-\ln v-1+\theta-\ln \theta-1)dx \le e_{0},$$
which together with Jensen's inequality yields that for any $t\in [0,T]$
\be\la{2.5}\al_1\le\int_N^{N+1}v(x,t)dx\le\al_2,\quad\al_1\le\int_N^{N+1}\theta(x,t)dx\le\al_2,
\ee
where $0<\al_1<\al_2$ are two roots of
\be z-\ln z-1=e_0.\nonumber\ee
Furthermore, multiplying \eqref{2.000} by $\frac{1}{Y_N(t)} $ and  integrating the resultant equality over $[N,N+1]$, we obtain after using \eqref{2.5}, \eqref{2.2}, and \eqref{2.13} that
 \be\ba\la{2.14}  \frac{1}{Y_N(t)}\int_N^{N+1}v(x,t) dx &\le C\int_N^{N+1}B_N(x,t) \left(1+\int_0^t\frac{\te+\frac{1}{2}v|\nb|^2}{B_N(x,\tau)Y_N(\tau)}d\tau\right) dx \\&
\le C +C\int_0^t\frac{1}{Y_N(\tau)}\int_N^{N+1}\left(\theta+\frac{1}{2}v|\nb|^2\right)dxd\tau
 \\&\le C+C\int_0^t\frac{1}{Y_N(\tau)}d\tau,\ea\ee
 which combined with  Gr\"onwall's inequality  and  \eqref{2.5} shows that for any $t\in [0,T],$
\be\la{2.15} \frac{1}{Y_N(t)} \le C.\ee
This together with    \eqref{2.000},  \eqref{2.13},  and  the fact that $C$ is independent of $N$ yields \eqref{cqq1} and thus finishes the proof of Lemma \ref{wq20}.  \thatsall

Now, with the similar arguments in \cite{py,cps}, we can obatin  the following   lower bound  of the temperature $\theta$.
\begin{lemma}\la{wq20x} There  exists a positive constant $C$ such  that  for any $(x,t)\in\O\times [0,T]$,
\be\ba\la{cqq1x}
\theta\ge C.
\ea\ee
\end{lemma}

\pf The proof is similar as those in \cite{py,cps}. We will sketch them here for completeness. Denoting by
  \be\ba\left(\theta >2\right)(t)=\left\{x\in\O|\theta(x,t)>2\right\}, ~~\left(\theta <1/2\right)(t)=\left\{x\in\O|\theta(x,t)<1/2 \right\},\nonumber\ea\ee
it follows from  \eqref{2.2} that
\be\ba e_0&\geq\int_{\left(\te<1/2\right)(t)}\left(\te-\ln \te-1\right)dx+\int_{\left(\te>2\right)(t)}\left(\te-\ln \te-1\right)dx\\&
\geq \left(\ln 2-1/2\right)\left|\left(\te<1/2\right)(t)\right|+\left(1-\ln 2\right)\left|\left(\te>2\right)(t)\right|\\
&\ge \left(\ln 2-1/2\right)  \left(\left|\left(\te<1/2\right)(t)\right|+ \left|\left(\te>2\right)(t)\right|\right)
, \nonumber\ea\ee
which shows that for any $t\in[0,T]$,
\be\ba\la{2.18}\left|\left(\te<1/2\right)(t)\right|+ \left|\left(\te>2\right)(t)\right| \le\frac{2e_0}{2\ln 2-1}. \ea\ee

Next, for any $p>2$, integrating \eqref{2.1}  multiplied  by $\te^{-2}\left(\te^{-1}-2\right)_{+}^{p}$ with $\left(\te^{-1}-2\right )_{+}\buildrel \Delta \over =\max\left\{\te^{-1}-2,0\right\}$ over $\O$, we obtain after using \eqref{cqq1} that
\bnn\ba &\frac{1}{p+1}\frac{d}{dt}\int_{\O}\left(\te^{-1}-2\right)_{+}^{p+1}dx+\int_{\O}\frac{ u_x^2+|\nw_x|^2+|\nb_x|^2}{v\te^2}\left(\te^{-1}-2\right)_{+}^pdx \\&\le\int_{\O}\frac{u_x}{v\te}\left(\te^{-1}-2\right)_{+}^p dx\\&
\le\frac{1}{2}\int_{\O}\frac{u_x^2}{v\te^2}\left(\te^{-1}-2\right)_{+}^pdx
+\frac{1}{2}\int_{\O}\frac{1}{ v}
\left(\te^{-1}-2\right)_{+}^p dx \\&
\le\frac{1}{2}\int_{\O}\frac{ u_x^2}{v\te^2}\left(\te^{-1}-2\right)_{+}^pdx+C\int_{\O}\left(\te^{-1}-2\right)_{+}^pdx,   \ea\enn
which along with   \eqref{2.18} leads to
\be\ba&\left\|\left(\te^{-1}-2\right)_{+}\right\|_{L^{p+1}(\O)}^{p} \frac{d}{dt}\left\|\left(\te^{-1}-2\right)_{+} \right\|_{L^{p+1}(\O)}
\le C\left\|\left(\te^{-1}-2\right)_{+}\right\|_{L^{p+1}(\O)}^{p} \nonumber  \ea\ee
with $C$ independent of $p$. This, in particular, implies that there exists some positive constant $C$ independent of $p$
such that
\bnn\sup_{0\le t\le T}\left\|\left(\te^{-1}-2\right)_{+}\right\|_{L^{p+1}(\O)}\le C.  \enn
 Letting $p\rightarrow +\infty$  and using \eqref{2.18}, it holds that
\be\ba\sup_{0\le t\le T}\left\|\left(\te^{-1}-2\right)_{+}\right\|_{L^{\infty}(\O)} \le C, \nonumber\ea\ee
which derives \eqref{cqq1x} and thus  finishes the proof of Lemma \ref{wq20x}. \thatsall

Now, we will prove the the following   upper bound of $v$, which is crucial for deducing the desired a priori estimates.

\begin{lemma} \la{upv} There exists a   positive constant $C $  such that   for any $(x,t)\in \Omega\times [0,T],$
\be \la{q2.49}v(x,t)\le C.\ee
Moreover, it holds that
\be \la{uiz1}\sup_{0\le t\le T}\int_{\Omega}(v^2_x+|\nb |^2)dx +\int_0^T\int_{\Omega}\left( (1+\te) v_x^2+u_x^2+|\nb_x|^2+|\nw_x|^2 \right)dxdt\le C.\ee
\end{lemma}

\pf
  First, multiplying \eqref{1.3} by $\nw$ and integrating the resultant
  equality over $\Omega  $ by parts, it holds that
\be\la{vd1}\ba  \frac{1}{2}\frac{d}{dt} \int_{\Omega} |\nw|^2dx +\int_{\Omega} \frac{|\nw_x|^2}{v}dx
 &\le C\int_{\Omega} |\nw_x|\left|\nb  \right|dx\\& \le C \xl(\int_{\Omega} \frac{|\nw_x|^2}{v}dx\xr)^{1/2}\xl(\int_{\Omega} v|\nb|^2dx\xr)^{1/2}\\& \le \frac12  \int_{\Omega} \frac{|\nw_x|^2}{v}dx+C,
\ea\ee
where in the last inequality one has used \eqref{2.2}.
Denoting by
\be \label{xiao1}  \ti V(t)\triangleq \int_{\Omega} \frac{|\nw_x|^2}{v}dx+V(t)+1\ee
where
 \be\label{xiao2} V(t) \triangleq \int_\O\left(\frac{u^2+|\nw|^2+v|\nb|^2}{2} \right)dx+\int_\O\left(\frac{ \theta_{x}^2}{v\theta^2 }+\frac{u_x^2+|\nw_x|^2+|\nb_x|^2}{v\theta}\right)dx,\ee
it deduces from \eqref{vd1} and \eqref{2.2} that
 \be\la{q2.45}\int_0^T\ti V(t) dt  \le C .\ee

Next, using \eqref{1.1}, we rewrite \eqref{1.2} as follows
\be\la{k1}  (\ln v)_{xt}=u_t+\left(\frac{\te}{v}\right)_x+\nb\cdot \nb_x.\ee
Adding  \eqref{k1}  multiplied by $(\ln v)_x$ and  \eqref{1.4}  multiplied   by $v\nb$ together,
 one obtains after  integrating the resultant equality over  $\Omega  $ by parts that
 \be\la{vg1}\ba &\frac12 \frac{d}{dt}\int_{\Omega} \left((\ln v)_x^2+v^2|\nb|^2\right)dx+\int_{\Omega} \left(  |\nb_x|^2+\frac{\te}{v}(\ln v)_x^2\right)dx\\& =\frac{d}{dt} \int_{\Omega} u(\ln v)_xdx +\int_{\Omega}\frac{u_x^2}{v}dx+\int_{\Omega} \frac{\te_x(\ln v)_x }{v}dx+\int_{\Omega} v\nw_x\cdot\nb dx.\ea\ee

Setting \be \label{xiao0} M_v(t)\triangleq 1+\max_{x\in \Omega}v(x,t),\ee
 the last term on the righthand side of \eqref{vg1} can be estiamted as follows
\be\la{vg12}\ba   \left|\int_{\Omega} v\nw_x \cdot \nb dx\right|
 \le CM_v(t)\int_{\Omega}   |\nw_x|| \nb | dx
\le  CM_v(t)\ti V(t) \ea\ee
due to   \eqref{vd1}.

In order to handle the third term on the righthand side of \eqref{vg1}, integration by   parts together with \eqref{2.1}  and \eqref{1.1}   gives
\be\la{vg100}\ba   \int_{\Omega} \frac{\te_x(\ln v)_x }{v}dx   &=-\int_{\Omega}  \ln v \left(\frac{\te_x }{v}\right)_xdx\\&=-\int_{\Omega}  \ln v \left(
(\theta-1)_{t}+  \frac\theta v
u_{x}-\frac{  u_{x}^{2}+|\nw_x|^2+|\nb_x|^2}{v}\right) dx\\&=-\left(\int_{\Omega} (\theta-1) \ln v dx\right)_{t} +\int_{\Omega} \frac{(\theta-1) u_{x}}{ v
} dx\\&\quad -\int_{\Omega} \frac{\theta u_{x}}{ v
}   \ln v  dx+\int_{\Omega} \frac{  u_{x}^{2}+|\nw_x|^2+|\nb_x|^2}{v}\ln v dx . \ea\ee
On the one hand, 
it holds
\be \la{q2.32} \ba    \int_{\Omega} \frac{\left(\te-1\right)}{v}u_xdx  &\le C \int_{\Omega} \frac{u_x^2}{v}dx+C \int_{\Omega} \left(\te-1\right)^2dx      \\
&\le C \int_{\Omega} \frac{u_x^2}{v}dx+ C\int_{(\te\le 2)(t)} \left(\te-1\right)^2dx+ C\int_{\O} \left(\te^{1/2}-2^{1/2}\right)_+^2\te dx      \\
&\le C \int_{\Omega} \frac{u_x^2}{v}dx+ C + C \max_{x\in\O} \left(\te^{1/2}-2^{1/2}\right)_+^2\int_{(\te > 2)(t)} \te dx      \\
 &\le  C\int_{\Omega} \frac{u_x^2}{v}dx +C+  C \left( \int_{(\te> 2)}|\te_x|\te^{-1/2}dx \right)^2   \\&\le  C\int_{\Omega} \frac{u_x^2}{v}dx + C+ C  \int_{\Omega} \frac{\te_x^2}{ \te^{2}v}dx  \int_{(\te>2)(t)} v\te dx   \\
&\le C\int_{\Omega} \frac{u_x^2}{v}dx +  C V(t)M_v(t)  +C,\ea\ee
where  one has  used \eqref{cqq1}  and the following estimates
\be\label{lls1} \int_{(\theta \le 2)(t)}(\theta-1)^2dx+\int_{(\theta>2)(t)} \theta dx \le C \int_\Omega(\theta-\ln \theta-1)dx \le C
\ee
owing to \eqref{2.2} and \eqref{2.18}.
On the other hand,
the straight calculations yield
\be \la{vg11} \ba  -\int_{\Omega} \frac{\te }{v}u_x\ln vdx  & = -\int_{\Omega} \frac{\te-1 }{v}u_x\ln vdx -\int_{\Omega} \frac{u_x\ln v }{v}dx    \\
&\le C\left(\int_{\Omega} \frac{u_x^2}{v}dx  +   \int_{\Omega} \left(\te-1\right)^2dx\right) \ln M_v(t )  \\
&\quad+C\int_{\Omega} \frac{u_x^2}{v}dx  +C\int_{\Omega} \ln^2 vdx\\
&\le C\left(\int_{\Omega} \frac{u_x^2}{v}dx  +     V(t)M_v(t)+1\right) \ln M_v(t),\ea\ee
where   one has used \eqref{cqq1}, \eqref{q2.32}, and the following fact
\be\ba \label{lls3} \int_{\Omega} \ln^2 vdx\le C\int_{(v\le 2)(t)} (v-1)^2dx +C \int_{(v> 2)(t)} vdx \le C\int_\Omega(v-\ln v-1)dx\le C\ea\ee
due to \eqref{2.2} and \eqref{cqq1}.

 Putting \eqref{q2.32} and \eqref{vg11} into \eqref{vg100} gives
\bnn\la{vg10}\ba   \int_{\Omega} \frac{\te_x(\ln v)_x }{v}dx    &\le
\frac{d}{dt} \int_\O -(\te-1)\ln vdx
 + C \left(\int_{\Omega} \frac{u_x^2+|\nb_x|^2}{v}dx   + \ti     V(t) M_v(t)\right)\ln M_v(t), \ea\enn which together with  \eqref{vg1}
  and  \eqref{vg12} leads to
  \be\la{vg13}\ba 
  &\frac12 \frac{d}{dt}\int_{\Omega} \left((\ln v)_x^2+|v\nb|^2\right)dx +\int_{\Omega} \left(  |\nb_x|^2+\frac{\te}{v}(\ln v)_x^2\right)dx \\
  &\le  \frac{d}{dt} B(t)  +C \left(\int_{\Omega} \frac{u_x^2+|\nb_x|^2}{v}dx +
   \ti V(t)M_v(t)\right)\ln M_v(t),\ea\ee
where
$$B(t) \triangleq \int_\O u(\ln v)_x dx-\int_\O (\te-1)\ln vdx $$
satisfies
\be\label{xiao3} \ba   B(t) &\le \frac{1}{8} \int_\O  (\ln v)_x^2 dx+C\int_\O u^2 dx+\int_\O |\te-1||\ln v|dx \\
&\le \frac{1}{8} \int_\O  (\ln v)_x^2 dx+C +\int_{(\te\le 2)(t)} (\te-1)^2 dx\\
&\quad+C\int_\O \ln^2vdx+\int_{(\te> 2)(t)}  \te  dx\ln M_v(t)\\
&\le \frac{1}{8} \int_\O  (\ln v)_x^2 dx+C +C\ln \int_\O  (\ln v)_x^2 dx\\
&\le \frac{1}{4} \int_\O  (\ln v)_x^2 dx+C\ea\ee
due to \eqref{2.2}, \eqref{lls1},  \eqref{lls3}, and the following fact
 \be\ba\label{xiao5} M_v(t)\le C+C\int_{\Omega} (\ln v)_x^2 dx .\ea\ee
Indeed, the direct calculations   combined with  \eqref{lls3} imply that
\bnn\ba  (v-2)_+ &\le C\left(\int_{(v>2)(t)}v^2 dx\right)^{1/2}\left(\int_{\Omega}\frac{v_x^2}{v^2}dx\right)^{1/2}
\\ &\le C \left(\int_{(v>2)(t)}v  dx  M_v(t)\right)^{1/2} \left(\int_{\Omega}(\ln v)^2_xdx\right)^{1/2}
\\ &\le C M^{1/2}_v(t) \left(\int_{\Omega}(\ln v)^2_xdx\right)^{1/2}, \ea\enn
which along with Young's inequality leads to \eqref{xiao5}.

Now, adding \eqref{1.2} multiplied by $u$ and \eqref{1.4} by $\nb$ together, integrating the resultant equality over  $(0,1)\times (0,t)$,   it follows from \eqref{cqq1}, \eqref{q2.32}, and \eqref{lls3} that
 \be\label{xxxiao6}\ba&\sup_{0\le s\le t}\int_{\Omega}  (u^2+v|\nb|^2)dx+\int_0^t \int_{\Omega} \frac{u_x^2+|\nb_x|^2}{v}dx ds\\
 &\le  C+C\xl|\int_0^t \int_{\Omega}\frac{\te u_x}{v} dxds\xr|\\
  &\le  C+C\xl|\int_0^t \int_{\Omega}\xl(\frac{(\te-1) u_x}{v}-\frac{(v-1) u_x}{v}+u_x\xr) dxds\xr|\\
  &\le  C+\frac{1}{2}\int_0^t \int_{\Omega} \frac{u_x^2}{v}dx ds+C \int_0^t \int_{\Omega} (\te-1)^2dxds+ C \int_0^t \int_{\Omega}(v-1)^2  dxds \\
  &\le C+\frac{1}{2}\int_0^t \int_{\Omega} \frac{u_x^2}{v}dx ds+C\int_0^t\ti V(s)M_v (s)ds\\
  &\quad+C \int_0^t \xl(\int_{(v\le 2)(t)}(v-1)^2  dx+M_v (s)\int_{(v> 2)(t)}v  dx\xr)ds\\
  &\le  C+\frac{1}{2}\int_0^t \int_{\Omega} \frac{u_x^2}{v}dx ds+C\int_0^t\ti V(s)M_v (s)ds.\ea\ee
This gives directly
\be\label{xxiao6}\int_0^t \int_{\Omega} \frac{u_x^2+|\nb_x|^2}{v}dx  ds\le C+C\int_0^t\ti V(s)M_v (s)ds,\ee
which combined with   \eqref{vg13} and \eqref{xiao3} yields
\be\la{2.45}\ba & \sup_{0\le s\le t}\int_{\Omega}  (\ln v)_x^2dx +\int_0^t \int_{\Omega} \left(  |\nb_x|^2+\frac{\te}{v}(\ln v)_x^2\right)dxds \\
&\le  C\ln \sup_{0\le s\le t}M_v(s) + C \int_0^t \ti V(s) M_v(s)ds \ln \sup_{0\le s\le t}M_v(s) +  C \\
&\le  C \left(2+ \int_0^t \ti V(s) M_v(s)ds  \right)\ln \left(2+ \int_0^t \ti V(s) M_v(s)ds  \right)  ,\ea\ee
where in the second inequality  one has used the following fact
\bnn fg\le  e^f-f-1+(1+g)\ln(1+g)-g,~~~~\mbox{ for any } f\ge 0,~g\ge 0, \enn
with
$$f=\frac12 \ln \sup_{0\le s\le t}M_v(s), \quad g= 2C \int_0^t \ti V(s) M_v(s)ds.$$

Then, the combination of  \eqref{xiao5} with \eqref{2.45} gives
 \bnn  M_v(t)  \le C\left(2+ \int_0^t \ti V(s) M_v(s)ds  \right)\ln \left(2+ \int_0^t \ti V(s) M_v(s)ds  \right),\enn
which along with   \eqref{q2.45}  and  Gr\"onwall's inequality derives \eqref{q2.49}.  Furthermore,  \eqref{uiz1} is  deduced directly form  \eqref{q2.49}, \eqref{2.45},  and \eqref{q2.45}.  The proof Lemma \ref{upv} is completed. $\hfill \Box$

\begin{lemma}\la{lemmy0} There  exists a positive constant $C$ such that
\be \la{aas1}\ba
&\sup_{0\le t\le T} \int_{\Omega}\left( |\nb_x|^2+|\nw_x|^2\right)dx \\&  +\int_0^T\int_{\Omega}\left(  |\nb_t|^2 +|\nb_{xx}|^2 +|\nw_t|^2+|\nw_{xx}|^2\right)dx dt\leq C.\ea\ee
\end{lemma}

\pf First,   rewriting  \eqref{1.3} as
\be\la{ppp}\ba \nw_t=\frac{\nw_{xx}}{v}-\frac{\nw_xv_x}{v^2}+\nb_x, \ea\ee multiplying \eqref{ppp} by $\nw_{xx},$ and integrating the resultant equality over $\Omega\times (0,T)$ by parts, we obtain after using
\eqref{uiz1}  and Cauchy  inequality that
\be\ba\la{bbb}&\frac{1}{2}\int_{\Omega}|\nw_x|^2dx+\int_0^T\int_{\Omega}\frac{|\nw_{xx}|^2}{v}dxdt   \\& \leq C+\frac{1}{4}\int_0^T\int_{\Omega}\frac{|\nw_{xx}|^2}{v}dxdt +C\int_0^T\int_{\Omega}\left(|\nb_x|^2+|\nw_x|^2v_x^2\right)dxdt \\&\leq C+\frac{1}{4}\int_0^T\int_{\Omega}\frac{|\nw_{xx}|^2}{v}dxdt +C\int_0^T\max_{x\in\Omega}|\nw_x|^2dt\\
&\le C+\frac{1}{4}\int_0^T\int_{\Omega}\frac{|\nw_{xx}|^2}{v}dxdt +C\int_0^T\xl(\int_\O |\nw_x|^2 dx\xr)^{1/2}\xl(\int_\O |\nw_{xx}|^2 dx\xr)^{1/2}dt\\
&\le C+\frac{1}{2}\int_0^T\int_{\Omega}\frac{|\nw_{xx}|^2}{v}dxdt, \ea\ee
where in the third inequality one has used the following fact
\be\ba\label{xiao6}  |f|^2 = \int_x^\infty (|f|^2)_ydy\le C\int_\O |f||f_x|  dx \le C \xl(\int_\O |f|^2 dx\xr)^{1/2}\xl(\int_\O |f_x|^2 dx\xr)^{1/2}\ea\ee
for any $f\in H^1(\O).$

 The combination of  \eqref{bbb}  with  \eqref{q2.49} leads to
\be\ba \la{ddd}\sup_{0\le t\le T}\int_{\Omega}|\nw_x|^2dx+\int_0^T\int_{\Omega}|\nw_{xx}|^2dxdt\leq C.\ea\ee
Thus, it follows from \eqref{ddd}, \eqref{ppp},   and \eqref{uiz1} that
\be\ba\la{ooo}\int_0^T\int_{\Omega}|\nw_t|^2dxdt &\leq
C\int_0^T\int_{\Omega} \left(|\nb_x|^2+|\nw_{xx}|^2 + v_x^2|\nw_x|^2\right)dxdt\\&\leq  C+ C\int_0^T\max_{x\in\Omega}|\nw_x|^2dt\\&\leq C .\ea\ee

Next,   rewriting \eqref{1.4} as
\be\la{kkk}\ba\nb_t=\frac{\nw_x}{v}+\frac{\nb_{xx}}{v^2}-\frac{\nb_x v_x}{v^3}-\frac{\nb u_x}{v},\ea\ee multiplying \eqref{kkk} by $\nb_{xx}$ and integrating the resultant equality  over $\Omega\times(0,T)$ by parts, we deduce from \eqref{uiz1}, \eqref{cqq1},  \eqref{2.2},   \eqref{xiao6},  and Cauchy  inequality that
\be\label{xiao4}\ba &\frac{1}{2}\int_{\Omega}|\nb_x|^2dx+\int_0^T\int_{\Omega}\frac{|\nb_{xx}|^2}{v^2}dxdt\\ &\leq C+\frac{1}{2}\int_0^T\int_{\Omega}\frac{|\nb_{xx}|^2}{v^2}dxdt +C\int_0^T\int_{\Omega}\left(|\nb_x|^2v_x^2+u_x^2|\nb|^2+|\nw_x|^2\right)dxdt \\
&\leq C+\frac{1}{2}\int_0^T\int_{\Omega}\frac{|\nb_{xx}|^2}{v^2}dxdt+C\int_0^T\max_{x\in\Omega}|\nb_x|^2dt +\max_{(x,t)\in\Omega\times[0,T]}|\nb|^2 \\
&\leq C+\frac{3}{4}\int_0^T\int_{\Omega}\frac{|\nb_{xx}|^2}{v^2}dxdt+C \int_0^T\int_{\Omega}|\nb_x|^2dxdt \\&\quad  +C\sup_{0\le t\le T}\int_{\Omega}|\nb|^2dx +\frac{1}{4}\sup_{0\le t\le T}\int_{\Omega}|\nb_x|^2dx \\&\leq C +\frac{3}{4}\int_0^T\int_{\Omega}\frac{|\nb_{xx}|^2}{v^2}dxdt+\frac{1}{4}\sup_{0\le t\le T}\int_{\Omega}|\nb_x|^2dx,\ea\ee
which along with \eqref{q2.49} implies
\be\ba\la{zzz}\sup_{0\le t\le T}\int_{\Omega}|\nb_x|^2dx+\int_0^T\int_{\Omega}|\nb_{xx}|^2dxdt\leq C.\ea\ee
This combined with \eqref{xiao6} and \eqref{uiz1} yields
\be\la{xxx} \max_{(x,t)\in \Omega\times [0,T]}|\nb|^2\le C \sup_{0\le t\le T}\xl(\int_\O |\nb|^2 dx\xr)^{1/2}\xl(\int_\O |\nb_x|^2 dx\xr)^{1/2}\le C.\ee

Finally, it follows from   \eqref{kkk}, \eqref{cqq1}, \eqref{zzz}, \eqref{uiz1},
\eqref{xxx}, and \eqref{xiao6} that
\bnn\ba\la{rrr}\int_0^T\int_{\Omega}|\nb_t|^2dxdt&\leq C\int_0^T\int_{\Omega}\left(|\nb_{xx}|^2+|\nb_x|^2v_x^2+|\nw_x|^2+|\nb|^2u_x^2\right)dxdt\\&\leq C+C\int_0^T \max_{x\in \Omega}|\nb_x|^2 \int_{\Omega}v_x^2dx  dt\\&\leq C+C\int_0^T\int_{\Omega}\left(|\nb_x|^2+|\nb_{xx}|^2\right)dxdt\\&\leq C.\ea\enn
   This together with \eqref{ddd}, \eqref{ooo},  and  \eqref{zzz}  gives \eqref{aas1} and thus  finishes the proof of Lemma \ref{lemmy0}. \thatsall

\begin{lemma}\la{nlemma70}
There  exists a positive constant $C$ such that
\be \la{7.1}\ba
&\sup_{0\le t\le T} \int_{\O} u_x^2 dx   +\int_0^T\int_{\O}\left( u_t^2 +u_{xx}^2 \right)dx dt\leq C.\ea\ee
\end{lemma}

\pf
First, we rewrite \eqref{1.2} as follows
\be\label{xiao8} u_t=\frac{u_{xx}}{v}-\frac{u_xv_x}{v^2}-\frac{\te_x}{v}+\frac{\te v_x}{v^2}-\nb\cdot\nb_x.\ee
Multiplying \eqref{xiao8} by $u_{xx}$ and integrating the result over $\O\times(0,T)$ by parts,  it holds
 \be\ba\la{7.2} &\frac{1}{2}\int_{\O}u_x^2dx+\int_0^T\int_{\O}\frac{ u_{xx}^2}{v}dxdt  \\&\leq C+\frac{1}{2}\int_0^T\int_{\O}\frac{ u_{xx}^2}{v}dxdt+C\int_0^T\int_{\O}\left(\te_x^2+\te^2v_x^2+|\nb|^2|\nb_x|^2+u_x^2v_x^2 \right)dxdt \\&\leq C+\frac{1}{2}\int_0^T\int_{\O}\frac{ u_{xx}^2}{v}dxdt+C\int_0^T\int_{\O}\te_x^2dxdt +C\int_0^T\sup_{x\in\O}\te^2\int_{\O}v_x^2dxdt \\&\quad
+C\sup_{(x,t)\in\O\times[0,T]}|\nb|^2\int_0^T\int_{\O}|\nb_x|^2dxdt
+C\int_0^T\sup_{x\in\O}u_x^2\int_{\O}v^2_xdxdt  \\
&\leq C+\frac{1}{2}\int_0^T\int_{\O}\frac{ u_{xx}^2}{v}dxdt+C\int_0^T\int_{\O}\te_x^2dxdt \\&\quad
+C\int_0^T\sup_{x\in\O}(\te-2)_{+}^2dt+C\int_0^T \left(\int_{\O}u_{x}^2dx\right)^{\frac{1}{2}}\left(\int_{\O}u_{xx}^2dx\right)^{\frac{1}{2}}dt  \\
&\le C+ \frac{3}{4} \int_0^T\int_{\O}\frac{ u_{xx}^2}{v}dxdt
+C_1\int_0^T\int_{\O}\te_x^2 dxdt,
 \ea\ee
where one has used  \eqref{uiz1}, \eqref{aas1},  \eqref{xxx}, \eqref{xiao6},  and the following  estimates:
\be\ba\label{xiao9}\sup_{x\in\O}(\te-2)_{+}^2&=\sup_{x\in\O}\left(\int_x^{\infty}\partial_y(\te-2)_{+} dy\right)^2
\le\left(\int_{(\te>2)(t)}|\te_y|dy\right)^2
\le C\int_{\O}\te_x^2dx  \ea\ee owing to \eqref{2.18}.

Then, motivated by \cite{L-L},   integrating \eqref{2.1} multiplied by $(\te-2)_{+}\buildrel\Delta\over=\max\limits_{x\in\O}\left\{\te-2,0\right\}$ over $\O\times(0,T)$, one obtains after using \eqref{aas1} and \eqref{cqq1} that
\be\ba\la{7.3}&\frac{1}{2}\int_{\O}(\te-2)_{+}^2dx+\int_0^T\int_{(\te>2)(t)}\frac{\te_x^2}{v}dxdt\\
&= \frac{1}{2}\int_{\O}(\te_{0}-2)_{+}^2dx+\int_0^T\int_{\O}\frac{u_x^2 + |\nw_x|^2 + |\nb_x^2|}{v}(\te-2)_{+}dxdt\\
&\quad-
\int_0^T\int_{\O}\frac{\te(\te-2)_{+}}{v}u_xdxdt\\&
\le C+C\int_0^T\sup_{x\in\O}\te\int_{\O}\left(|\nb_x|^2+|\nw_x|^2+(\te-2)_{+}^2+u_x^2\right)dxdt\\&
\le C+C\int_0^T\sup_{x\in\O}\te\left(\int_{\O}(\te-2)_{+}^2dx+\int_{\O}u_x^2dx+1\right)dt.
\ea\ee
Furthermore, it follows from \eqref{2.2} and \eqref{q2.49} that
\be\ba\int_0^T\int_{\O}\te_x^2dxdt&=\int_0^T\int_{(\te>2)(t)}\te_x^2dxdt+\int_0^T\int_{(\te\le 2)(t)}\te_x^2dxdt\\&
\le C\int_0^T\int_{(\te>2)(t)}\frac{\te_x^2}{v} dxdt+ C\int_0^T\int_{(\te\le 2)(t)}\frac{\te_x^2}{v\te^2} dxdt\\&
\le C\int_0^T\int_{(\te>2)(t)}\frac{\te_x^2}{v} dxdt+ C,\nonumber\ea\ee
which together with \eqref{7.3} yields
 \be\ba\la{7.4}&\frac{1}{2}\int_{\O}(\te-2)_{+}^2dx+C_2\int_0^T\int_{\O} \te_x^2 dxdt\\&
\le  C+C\int_0^T\sup_{x\in\O}\te\left(\int_{\O}(\te-2)_{+}^2dx+\int_{\O}u_x^2dx+1\right)dt.
\ea\ee

Thus, one obtains after adding \eqref{7.4} multiplied by $2C_2^{-1}C_1$ to \eqref{7.2} that
\be\ba\label{xiao7}  &\int_{\O}\left(u_x^2+(\te-2)_{+}^2\right)dx+\int_0^T\int_{\O}\left(u_{xx}^2+\te_x^2\right)dxdt\\
&\le C+C\int_0^T\sup_{x\in\O}\te\left(\int_{\O}(\te-2)_{+}^2dx+\int_{\O}u_x^2dx+1\right)dt.   \ea\ee
The straight calculations together with \eqref{q2.49} and \eqref{lls1} imply that
\be\la{jdte1}\ba\nonumber
\te^{1/2} \leq \left(\te^{1/2}-2^{1/2}\right)_{+}+C &\leq \int_{(\te \geq 2)(t)} \frac{|\te_x|}{\te^{\frac{1}{2}}} dx + C \\&\leq C\left(\int_{(\te \geq 2)(t)} \frac{\te_x^2}{v\te^2} dx\right)^{1/2} \left(\int_{(\te \geq 2)(t)} v\te dx\right)^{1/2}+C \\&\leq C\left(\int_{(\te \geq 2)(t)} \frac{\te_x^2}{v\te^2} dx\right)^{1/2}+C,
\ea\ee
which along with \eqref{2.2} yields
\be\la{wdi1}\ba
\int_{0}^{T} \sup_{x\in\O} \te dt \leq  C+C\int_0^T \int_{\O} \frac{\te_x^2}{v\te^2} dx dt \leq C.
\ea\ee
The Gr\"onwall's inequality together with \eqref{xiao7} and  \eqref{wdi1} leads to
\be\ba\la{7.5} \sup_{0\le t\le T}\int_{\O}\left(u_x^2+(\te-2)_{+}^2\right)dx+\int_0^T\int_{\O}\left(u_{xx}^2+\te_x^2\right)dxdt\le C.  \ea\ee

Finally,  it is easy to derive from \eqref{xiao8}, \eqref{cqq1}, \eqref{7.5}, \eqref{uiz1}, \eqref{xxx}, \eqref{xiao6}, \eqref{xiao9}, and \eqref{aas1} that
\bnn\ba
\int_0^T\int_{\O}u_t^2dxdt&\leq C\int_0^T\int_{\O}\left(u_{xx}^2+u_x^2v_x^2+\te_x^2+\te^2v_x^2+|\nb|^2|\nb_x|^2\right)dxdt\\
&\leq C+C\int_0^T\sup_{x\in\O}\left(u_x^2  +(\te-2)_+^2  \right)\int_{\O} v_x^2 dxdt\\
&\le C+C\int_0^T \left(\int_{\O}u_{x}^2dx+\int_{\O}u_{xx}^2dx++\int_{\O}\te_{x}^2dx\right)dt \\&\le C,
\ea\enn
which combined  with  \eqref{7.5}   gives \eqref{7.1} and thus completes the proof of \lemref{nlemma70}.\thatsall

\begin{lemma}\la{nlemma80}There exists a positive constant $C$ such that \be\ba\la{eq1} \sup_{(x,t)\in\O\times[0,T]}\te(x,t)+\sup_{0\le t\le T}\int_{\O } \te_x^2dx+\int_0^T\int_{\O}  \left( \te_t^2+\te_{xx}^2\right)dxdt\le C .  \ea\ee
\end{lemma}

\pf First, multiplying \eqref{2.1} by
$ \te_t$ and integrating the resultant equality over $\O$, it holds
\be\ba\la{8.2} & \int_{\O}  \te_t^2dx+ \frac{1}{2} \left(\int_{\O} \frac{\theta_{x}^2}{v}dx\right)_t\\&=- \frac{1}{2}\int_{\O}\frac{\theta_{x}^2u_x}{v^2}dx+\int_{\O} \frac{ \theta_t\left(-\te u_x+ u_x^2+|\nw_x|^2+|\nb_x|^2\right)}{v }dx\\
&\le C\sup_{x\in\O} (|u_x|)\int_{\O} \te_x^2dx +\frac{1}{2}\int_{\O} \te_t^2dx+C\int_{\O}\te^{2}u_x^2dx\\&\quad+C\int_{\O}\left(u_x^4+|\nw_x|^4+|\nb_x|^4\right)dx\\
&\le C \left(\int_{\O}\te_x^2dx\right)^2+\frac{1}{2}\int_{\O}\te_t^2dx\\
&\quad+C\sup_{x\in\O} \xl(  (\te-2)_+^{2}+  u_x^2 +|\nw_x|^2+|\nb_x|^2 \xr) \xl(\int_{\O}\left(u_x^2+|\nw_x|^2+|\nb_x|^2\right)dx+1\xr)\\
&\le  C \left(\int_{\O}\te_x^2dx\right)^2+\frac{1}{2}\int_{\O}\te_t^2dx
+C\int_{\O}\left(\te_x^2+u_{xx}^2+|\nw_{xx}|^2+|\nb_{xx}|^2\right)dx+C ,\ea\ee
where one has used  \eqref{cqq1}, \eqref{xiao9}, \eqref{xiao6},
\eqref{aas1},   and \eqref{7.5}.
 Thus,   Gr\"onwall's inequality together with \eqref{8.2}, \eqref{7.5}, and \eqref{aas1} implies  that
\be\ba\la{8.5} \sup_{0 \le t\le T}\int_{\O} \theta_{x}^2 dx+\int_0^T\int_{\O} \te_t^2dxdt\le C, \ea\ee
which combined with \eqref{xiao9}  gives
\be\la{8.6}\sup_{(x,t)\in\O\times[0,T]}\te(x,t)\le C.\ee

Finally, it follows from \eqref{2.1} that
\bnn\ba \frac{\te_{xx}}{v}= \frac{\te_x v_x}{v^2}- \frac{ u_x^2+|\nb_x|^2+|\nw_x|^2}{v}+ \frac{ \te u_x}{v}+\te_t,\ea\enn
which together with \eqref{q2.49},  \eqref{cqq1},  \eqref{8.6},  \eqref{uiz1}, \eqref{aas1}, \eqref{7.1}, \eqref{xiao6},  and   \eqref{8.5} yields
\bnn\ba\la{8.8}\int_0^T\int_{\O} \te_{xx}^2dxdt  & \le C\int_0^T\int_{\O} \left(\te_x^2v_x^2+u_x^4+|\nb_x|^4+|\nw_x|^4+u_x^2+\te_t^2\right)dxdt\\&
\le C+ C\int_0^T\sup_{x\in\O}\te_x^2 dt+\int_0^T\sup_{x\in\O}\left(|u_x|^2+ |\nb_x|^2+ |\nw_x|^2\right)dt \\&
\le C+\frac12\int_0^T \int_{\O} \te_{xx}^2 dx   dt.\ea\enn
 Combining this with \eqref{8.5}--\eqref{8.6}  proves \eqref{eq1}
and  finishes the proof of \lemref{nlemma80}.\thatsall

 \end{document}